\documentclass{amsart}
\usepackage{amsmath}
\usepackage{amsfonts}
\usepackage{amsthm}
\usepackage{amssymb}
\usepackage[english]{babel}
\usepackage[all]{xy}
\usepackage{color}

\theoremstyle{definition}

\newcommand{\dC}{{\mathcal C}}

\newcommand{\dE}{{\mathcal E}}
\newcommand{\dF}{{\mathcal F}}

\newcommand{\dO}{{\mathcal O}}

\newcommand{\dT}{{\mathcal T}}

\newtheorem{defi}{Definition}[section]
\newtheorem{remark}[defi]{Remark}
\newtheorem{notation}[defi]{Notation}

\newtheorem{definition}[defi]{Definition}

\theoremstyle{plain}

\newtheorem{theorem}[defi]{Theorem}
\newtheorem{corollary}[defi]{Corollary}
\newtheorem{lemma}[defi]{Lemma}
\newtheorem{proposition}[defi]{Proposition}

\input xy
\xyoption{all}

\title{Models of torsors and the fundamental group scheme}

\author{Marco Antei}\address{Marco Antei, Laboratoire J.A.Dieudonn\'e
UMR CNRS-UNS N${}^o$7351
Universit\'e de Nice Sophia-Antipolis
Parc Valrose
06108 NICE Cedex 2}\thanks{Marco Antei remercie le projet TOFIGROU (ANR-13-PDOC-0015-01).} 
\author{Michel Emsalem}\thanks{Michel Emsalem a re\c{c}u le soutien du Labex CEMPI (ANR-11-LABX-01)} \address{Michel Emsalem, Laboratoire Paul Painlev\'e, U.F.R. de Math\'ematiques, Universit\'e des Sciences et des Technologies de Lille 1, 59 655 Villeneuve d'Ascq, France}

\date{}
\begin{document}

\begin{abstract}
Given a relative faithfully flat pointed scheme over the spectrum of a discrete valuation ring $X \to S$ this paper is motivated by the study of the natural morphism from the fundamental group scheme of the generic fiber $X_\eta $ to the generic fiber of the fundamental group scheme of $X$. Given a torsor $T \to X_\eta $ under an affine group scheme $G$ over the generic fiber of $X$, we address the question to find a model of this torsor over $X$, focusing in particular on the case where $G$ is finite. We obtain partial answers to this question, showing for instance that, when $X$ is integral and regular of relative dimension $1$, such a model exists on some model of $X_{\eta}$ obtained by performing a finite number of N\'eron blow-ups along a closed subset of the special fiber of $X$. In the first part we show that the relative fundamental group scheme of $X$ has an interpretation as the Tannaka Galois group of a tannakian category constructed starting from the universal torsor. 
\end{abstract}  
\maketitle

\textbf{Mathematics Subject Classification. Primary: 14L30, 14L15. Secondary: 11G99.}\\\indent
\textbf{Key words:} torsors, quasi-finite group schemes, fundamental group scheme.

\tableofcontents
\bigskip
\indent \textbf{Acknowledgements} 
We would like to thank Hélène Esnault for very useful and helpful discussions and for pointing out a problem in a preliminary version. We also would like to thank Carlo Gasbarri for very interesting conversation on this subject. 
\section{Introduction}

\subsection{Aim and scope}\label{sez:Aim}
Let $S$ be a Dedekind scheme of dimension one and $\eta=Spec(K)$ its generic point; let $X$ be a scheme, $f:X\to S$ a faithfully flat morphism of finite type and $f_{\eta}:X_{\eta}\to \eta$ its generic fiber. Assume we are given a finite $K$-group scheme $G$ and a $G$-torsor $Y\to X_{\eta}$. So far the problem of extending the $G$-torsor $Y\to X_{\eta}$ has consisted in finding a finite and flat  $S$-group scheme $G'$ whose generic fibre is isomorphic to $G$ and  a $G'$-torsor $T\to X$ whose generic fibre is isomorphic to $Y\to X_{\eta}$ as a $G$-torsor.  Some solutions to this problem, from Grothendieck's first ideas until nowadays, are known in some particular relevant cases that we briefly recall:  Grothendieck proved that, possibly after extending scalars, the problem has a solution when $G$ is a constant finite group, $S$ is the spectrum of a complete discrete valuation ring with algebraically closed residue field of positive characteristic $p$, with $X$ proper and smooth over $S$ with geometrically connected fibers  and $p\nmid |G|$ (\cite{SGA1}, Expos\'e X);  when $S$ is the spectrum of a discrete valuation ring of residue characteristic $p$, $X$ is a proper and smooth curve over $S$ then Raynaud suggested a solution, possibly after extending scalars, for $|G|=p$  (\cite{Ray3}, \S 3); when $S$ is the spectrum of a  discrete valuation ring $R$ of mixed characteristic $(0,p)$ Tossici provided a solution, possibly after  extending scalars, for $G$ commutative when $X$ is a regular scheme, faithfully flat over $S$, with some extra assumptions on $X$ and $Y$  (\cite{Tos}, Corollary 4.2.8). 
Finally in \cite{Antei2}, \S 3.2 and \S 3.3  the first author provided a solution for $G$ commutative, when  $S$ is a connected Dedekind scheme and $f:X\to S$ is a smooth morphism  satisfying additional assumptions (in this last case it is not needed to extend scalars) and in \cite{Antei4} the case where $G$ is solvable is treated.  However a general solution does not exist. Moreover we know that it can even happen that $G$ does not  admit a finite and flat model (see \cite{Mi2}, Appendix B, Proposition B.2 for the positive equal characteristic case or \cite{Ray}, \S 3.4 for the mixed characteristic case). What is always true is that $G$ admits at least an affine, quasi-finite (then of finite type, according to our conventions, see \S \ref{sez:Conve}), flat $R$-group scheme model. Indeed $G$  is isomorphic to a closed subgroup scheme of some $GL_{n,K}$ (\cite{WW}, \S 3.4) then it is sufficient to consider its schematic closure  in $GL_{n,S}$. In this paper we explain how to solve the problem of extending any $G$-torsor when $G$ is any algebraic group scheme over $K$ up to a modification of $X$. When $X$ is a relative curve this \emph{modification} is represented by a finite number of N\'eron blow-ups of $X$ in a closed subscheme of the special fiber of $X$. For more precise statements we refer the reader to \S \ref{sez:tutta}. The most interesting case is certainly the case where $G$ is finite. If we were able to prove that every finite and pointed torsor over $X_{\eta}$ admits a model over $X$ then the natural morphism $\varphi:\pi(X_{\eta},x_{\eta})\to \pi^{\rm qf}(X,x)_{\eta}$ (which is always faithfully flat, \cite[\S 7]{AEG}) between the fundamental group scheme of $X_{\eta}$ to the generic fiber of the quasi-finite fundamental group scheme of $X$ would be an isomorphism. It is known that $\varphi$ becomes an isomorphism when we restrict to the abelianized fundamental group scheme (cf. \cite{Antei2}). Here we find a partial answer, extending \emph{all} finite torsors, but instead of providing a model over $X$ we provide a model over some $X'$ obtained slightly modifying $X$, as explained. In order to approach the question from a different point of view it would be of great interest to have a tannakian description for $\pi^{\rm qf}(X,x)$, $\pi(X,x)$  and their universal torsors $\widehat{X}^{\rm qf}\to X$ and $\widehat{X}\to X$. In \cite{MS} Mehta and Subramanian provided a first construction which works only for schemes defined over some non-noetherian Pr\"ufer rings whose function field is algebraically closed. In \S \ref{sez:Tanna} we give a different tannakian description simply choosing the category of vector bundles on $X$  trivialized by the universal torsor, whose existence is now known. An intrinsic description, independent from the existence of the universal torsor, would be strongly appreciated.     

\subsection{Notations and conventions}
\label{sez:Conve}  Let $S$ be any scheme, $X$ a $S$-scheme, $G$ a (faithfully) flat, affine\footnote{We do not need $G$ to be affine in order to define a $G$-torsor, but it is the only type we encounter in this paper.} $S$-group scheme and $Y$ a $S$-scheme endowed with a right action $\sigma:Y\times G\to Y$. A $S$-morphism $p : Y  \to X$ is said to
be a $G$-torsor if it is 
faithfully flat, $G$-invariant and the canonical morphism $(\sigma,pr_Y):Y\times G\to Y\times_X Y$ is an isomorphism. Let $H$ be a flat, affine $S$-group scheme and $q:Z\to X$  a $H$-torsor; a morphism  between two such torsors is a pair $(\beta,\alpha):(Z,H)\to(Y,G)$ where $\alpha:H\to G$  is a $S$-morphism of group schemes, and $\beta:Z\to Y$ is a $X$-morphism of schemes such that  the following diagram commutes

$$\xymatrix{Z\times H \ar[r]^{\beta\times \alpha} \ar[d]_{H\text{-}action} & Y\times G\ar[d]^{G\text{-}action}\\ Z\ar[r]_{\beta} & Y}$$ (thus $Y$ is isomorphic to the contracted product $Z\times^H G$  through $\alpha$, cf. \cite{DG}, III, \S 4, 3.2). In this case we say that $Z$ precedes $Y$. 

When $S$ is irreducible, $\eta$ will denote its generic point and $K$ its function field $k(\eta)$. Any $S$-scheme whose generic fiber is isomorphic to some $K$-scheme $T_{\eta}$ will be called a model of $T_{\eta}$. Any morphism which is generically an isomorphism will be called a model map. 

Throughout the whole  paper a  morphism of schemes $f:Y\to X$ will be said to be quasi-finite if it is of finite type and for every point $x\in X$ the fiber $Y_x$ is a finite set. Let $S$ be any scheme and $G$ an affine  $S$-group scheme. Then we say that $G$ is  a finite (resp. quasi-finite) $S$-group scheme if the structural morphism $G\to S$ is finite (resp. quasi-finite).

Let $R$ be a commutative ring with unity and $G$ an affine and flat $R$-group scheme; we denote by $Rep_{R,tf} (G)$ the category of finitely generated $R$-linear representations of $G$ and by $Rep_{R,tf}^0 (G)$ the full subcategory of $Rep_{R,tf} (G)$ whose objects are free (as $R$-modules).

\section{A tannakian construction}
\label{sez:Tanna}

Throughout this section $R$ will always be a discrete valuation ring, with uniformising element $\pi$, with field of fractions $K:=Frac(R)$ and residue field $k:=R/\pi R$. The generic and special points of $Spec(R)$ will often be denoted by $\eta$ and $s$ respectively.

\begin{lemma}\label{lemLimit} Let $X$ be a noetherian scheme over $R$, $T = \varprojlim_{i\in I} T_i$ a projective limit of $X$-schemes $f_i : T_i \to X$ affine over $X$. We assume that for all $i\in I$, $H^0 (T_i, \dO _{T_i}) = R$. Then $H^0 (T, \dO _T)=R$. 
\end{lemma}
\proof  This follows from \cite{H77}, III, Poposition 2.9.
\endproof

\begin{lemma}\label{lemESSAI} Let $j : T\to Spec(R)$ be a surjective faithfully flat morphism and let us assume that the generic fiber  $T_{\eta}$ of $T$ is such that  $H^0 (T_{\eta}, \dO _{T_{\eta}})=K$. Then $H^0 (T, \dO _T)=R$.
\end{lemma}
\proof

We first observe that either $H^0 (T, \dO _T)=R$ or $H^0 (T, \dO _T)=K$. Indeed 
 $$R \subset H^0 (T, \dO _T) \subset H^0 (T, \dO _T)\otimes _R K \simeq  H^0 (T_\eta , \dO _{T_\eta })=K,$$
the last isomorphism being a consequence of \cite[III, Proposition 9.3]{H77}, 
and it is not difficult to check that a $R$-algebra containing $R$ and contained in  $K$ is either $R$ or $K$. However if $H^0 (T ,\dO _T)=K$, then the canonical factorisation of  $f: T \to Spec (R)$ into $T \to Spec (\dO_T(T)) \to Spec (R)$ would not be surjective.
\endproof

We apply these remarks to the theory of the (quasi-finite) fundamental group scheme and its universal torsor that we briefly recall: 

\begin{definition}
Let $X\to S$ be a morphism of schemes endowed with a section $x\in X(S)$. We say that $X$ has a fundamental group scheme (resp. a quasi-finite fundamental group scheme) if there exists a $S$-group scheme $\pi(X,x)$ (resp. $\pi^{\rm qf}(X,x)$) and a pointed $\pi(X,x)$-torsor $\widehat{X}$ (resp. $\pi^{\rm qf}(X,x)$-torsor $\widehat{X}^{\rm qf}$) such that for any finite (resp. quasi-finite) torsor $Y\to X$ over $X$, pointed over $x$, there is a unique morphism of pointed torsors $\widehat{X}\to Y$ (resp. $\widehat{X}^{\rm qf}\to Y$).  
\end{definition}

In  \cite{AEG}, \S 4,\S 5.1 and  \S 5.2 we proved the following existence theorems:
 
 \begin{theorem}\label{teoCofilDed}
Let $S$ be a Dedekind scheme, $X\to S$ a faithfully flat morphism locally of finite type and $x\in X(S)$ a section. Let us moreover assume that one of the following assumptions is satisfied:
\begin{enumerate}

\item for every $s\in S$ the fiber $X_s$ is reduced; 
\item for every $z\in X\backslash X_{\eta}$ the local ring $\mathcal{O}_{X,z}$ is integrally closed.
\end{enumerate}
Then $X$ has a fundamental group scheme. 
\end{theorem}

\begin{theorem}\label{teoCofilDed2}
Let $S$ be a Dedekind scheme, $X\to S$ a faithfully flat morphism locally of finite type and $x\in X(S)$ a section. Let us moreover assume that $X$ is integral and normal and that for each $s\in S$ the fiber $X_s$ is normal and integral. Then $X$ has a quasi-finite fundamental group scheme. 
\end{theorem}

\begin{notation} In order to simplify the exposition from now on we will only consider the case where $X$ satisfies the assumptions of Theorem \ref{teoCofilDed2}. However the statements proved for the universal $\pi^{\rm qf}(X,x)$-torsor will also hold for the universal $\pi(X,x)$-torsor and the proofs are exactly the same. 
\end{notation}

 \begin{definition}We say that a quasi-finite $G$-torsor over $X$ pointed over $x$ is quasi-Galois if the canonical morphism $\rho:\pi(X,x)^{\rm qf}\to G$ is schematically dominant (i.e. $R[G]\to R[\pi(X,x)^{\rm qf}]$ is injective). It will be furthermore called Galois, if $\rho:\pi(X,x)^{\rm qf}\to G$ is faithfully flat.
 \end{definition}
 
 \begin{lemma}\label{lemQGalois} Let us assume that $H^0(X,\dO_X)=R$. Let $G$ be a quasi-finite and flat group scheme and $T \to X$ a quasi-Galois $G$-torsor pointed over $x$. Then $H^0(T,\dO_T)=R$.
 \end{lemma}

 \proof It is sufficient to notice that the generic fiber $T_{\eta}$ of $T$ is Galois over $X_{\eta}$ (cf. \cite[\S 7]{AEG}), hence  $H^0(T_{\eta},\dO_{T_{\eta}})=K$ (cf. \cite[Chapter 2, Proposition 3]{Nor2}). The conclusion follows by Lemma \ref{lemESSAI}. \endproof
 
 \begin{corollary}\label{corFonda} Let us assume that $H^0(X,\dO_X)=R$. Then $ H^0 (\widehat{ X}^{\rm qf}, \dO_{\widehat{ X}^{\rm qf}})=R$.
 \end{corollary}
 \proof First we observe that for any quasi-finite and flat $R$-group scheme $G$, any $G$-torsor over $X$ pointed  over $x$ is preceded by a quasi-Galois torsor: it is sufficient to factor the canonical morphism $\pi^{\rm qf}(X,x)\to G$ as $\pi^{\rm qf}(X,x)\to G'\to G$ where $\pi^{\rm qf}(X,x)\to G'$ is a schematically dominant morphism and $G'\to G$ is a closed immersion. Then the contracted product $\widehat{X}^{\rm qf}\times^{\pi^{\rm qf}(X,x)} G'$ gives the desired quasi-Galois torsor. Hence the universal torsor is isomorphic to a projective limit of quasi-Galois torsor and the conclusion follows using Lemma \ref{lemQGalois}, 
 and Lemma \ref{lemLimit}. 
 \endproof 
 
 \begin{theorem}\label{theoTann} Assumptions being as in Corollary \ref{corFonda}, the universal $\pi(X,x)^{\rm qf}$-torsor $\widehat{X}^{\rm qf}\to X$ induces an equivalence of categories $ Rep_{R,tf}^0 (\pi^{\rm qf} (X,x)) \to \dT^{\rm qf}$ where $\dT^{\rm qf} $ denotes the category of vector bundles on $X$ trivialized by $\widehat{X}^{\rm qf}\to X$. 
\end{theorem}

Let $\theta :X\to Spec(R)$ denote the structural morphism and let $j:T\to X$ be any $G$-torsor for some affine and flat $R$-group scheme $G$, then we denote by $G-Vect_T$ the category of $G$-sheaves over $T$ whose objects are locally free as $\dO_T$-modules; it is known that the functor $j^\ast :Vect _X\to G-Vect _T $ is an equivalence of categories and we denote by $\rho :G-Vect _T \to Vect _X$ a quasi-inverse. As usual, we naturally associate to $j:T\to X$ the functor $F_T :Rep_{R,tf} ^0 G \to \dT$, where $\dT $ denotes the category of vector bundles on $X$ trivialized by $T\to X$, given by $\rho \circ (j^\ast \circ \theta ^\ast )$. Thus Theorem \ref{theoTann} is a consequence of the following more general statement which somehow generalizes an analog result for torsors over fields (cf. \cite[Chapter II, Proposition 3]{Nor2}):
 
 \begin{lemma}\label{lemFully} Let $G\to Spec (R)$ be a flat affine group scheme and $j:T \to X$ a $G$-torsor such that $H^0 (T, \dO_T)=R$. The functor $F_T : Rep_{R,tf} ^0  G \to \dT$  is an equivalence of categories. 
 \end{lemma}
 
 \proof  First we observe that the functor $F_T$ is fully faithful if and only if  $ (j^\ast \circ \theta ^\ast )$ is faithfully flat:
let $V_1, V_2$ be two objects of $Rep_{R,tf} ^0  G$. Then $Hom_{(Rep_{R,tf} ^0  G)} (V_1,V_2) = (V_1^{\vee } \otimes _R V_2)^G$. Analogously if  $\dF_1,\dF_2$ are two objects of  $G-Vect _T$, then $Hom _{G-Vec _T}(\dF_1,\dF_2)= H^0 (T,\dF_1^{\vee}\otimes _{\dO _T}\dF_2)^G$. Thus $F_T$ is fully faithful if and only if, for any object $W$ of $Rep_{R,tf} ^0  G$, the natural map $$(\dagger)\qquad W^G \to H^0 (T,j^*\theta^* (W))^G$$ is an isomorphism. We have the following sequence of isomorphisms (by means of the projection formula):

$$ H^0 (T,j^*\theta^* (W)) \simeq H^0 (Spec(R),(\theta \circ j)_\ast (\theta \circ j)^*(W))\simeq H^0 (Spec(R),(\theta \circ j)_\ast \dO_T \otimes_{R} W)\simeq $$
$$\simeq  H^0 (Spec(R),(\theta \circ j)_\ast \dO_T) \otimes_{R} W\simeq  H^0 (T, \dO _T) \otimes _R W $$
as representations of $G$, then $(H^0 (T,j^*\theta^* (W)))^G \simeq (H^0 (T, \dO _T) \otimes _R W)^G=W^G$ since we assumed $H^0 (T, \dO _T)=R$, whence the desired isomorphism $(\dagger)$. 
In order to prove the essential surjectivity, we argue as follows: let us take a vector bundle $E$ on $X$ trivialized by $j:T \to X$. That implies the existence of a finitely generated free $R$-module $M$ such that  $\dE:= j^\ast E \simeq (\theta \circ j)^\ast M$. Again applying the projection formula we obtain

$$ (\theta \circ j)_\ast (\theta \circ j)^\ast M= (\theta \circ j)_\ast \dO_T\otimes_R  M=M.$$
It follows that $\dE \simeq (\theta \circ j)^\ast (\theta \circ j)_\ast (\theta \circ j)^\ast M \simeq (\theta \circ j)^\ast (\theta \circ j)_\ast \dE \simeq (\theta \circ j)^\ast H^0 (T, \dE )$. We now observe that the previous isomorphism $(\theta \circ j)^\ast (\theta \circ j)_\ast \dE \to \dE $  is $G$-equivariant  and thus $F_T(H^0 (T,  \dE))\simeq E$.
 \endproof

\begin{remark}Lemma \ref{lemFully} can be generalized further as follows: let $R$ be any commutative and unitary ring, $q:T\to Spec(R)$ a morphism of scheme such that $H^0 (T, \dO _T)=R$. Let moreover $G$ be any flat and affine $R$-group scheme, acting on $T$ and let $\dF$ be any $G$-sheaf, trivial as a $\dO _T$-module. Then $H^0 (T, \dF )$ is a $R$-linear representation of $G$ and $\dF \simeq H^0 (T, \dF ) \otimes _R \dO _T$ as $G$-sheaves.
\end{remark}

\begin{remark}Notations being as in Lemma \ref{lemFully}, the tannakian category $\dC$ over $R$ (cf. \cite{DH} for a modern and detailed overview) associated to $\dT$ is the category of those $\dO_X$-modules whose pullback over $T$ is isomorphic, as $\dO _T$-module, to a finite direct sum of $\dO _T$ and $\dO _T /\pi ^n$, where $\pi$ denotes a uniformizer of $R$ and $n$ a natural integer; in this way  $\dT$ would coincide with the full subcategory $\dC ^0$ of $\dC$ of rigid objects of $\dC$, i.e. objects isomorphic to their double dual. It would be very interesting and useful to have an inner description of the objects of $\dT$ (or equivalently of $\dC$) independent from the universal torsor. 
\end{remark}

\section{Existence of a model}
\label{sez:tutta}

\subsection{Quotients and N\'eron blow-ups}
\label{sez:almost}

In this section we are going to recall some results ensuring the existence of quotients of schemes under the action of some group schemes, under certain assumptions. Those results are essentially contained in \cite{SGA1}, Expos\'e V, Th\'eor\`eme 7.1 and \cite{Ray0}, Th\'eor\`eme 1, (iv), for the finite case and \cite{Ray0}, Th\'eor\`eme 1, (v) and \cite{ANA}, Th\'eor\`eme 7, Appendice 1, for the quasi-fini case. The fact that quotients (under the action of finite type group schemes) commute with base change is ensured by \cite[Expos\'e IV, 3.4.3.1]{SGA1}.

\begin{theorem}\label{teoQuozReno}
Let $T$ be a locally noetherian scheme, $Z$ a $T$-scheme locally of finite type, $H$ a flat $T$-group scheme acting on $Z$ in such a way that $Z\times_T H\to Z\times_T Z$ is a closed immersion. Then if one of the following conditions is verified 
\begin{enumerate}
\item $H\to T$ is finite and $Z\to T$ is quasi-projective,
\item $H\to T$ is quasi-finite and $Z\to T$ is quasi-finite, 
\end{enumerate}
the sheaf $(Z/H)_{fpqc}$ is represented by a scheme $Z/H$.  Moreover the canonical morphism $Z\to Z/H$ is faithfully flat and the natural morphism $Z \times _T H \to Z\times _{Z/H} Z $ is an isomorphism. 
\end{theorem}

\begin{theorem}\label{teoQuozReno2}
Let $T$ be any  locally noetherian scheme, $Z$ a  $T$-scheme locally of finite type, $H$ a flat $T$-group scheme acting on $Z$ such that the natural morphism $Z\times_T H\to Z\times_T Z$ is a closed immersion. Then there exists a largest open $U$ of $Z$ for which the sheaf $(U/H)_{fpqc}$ is represented by a scheme $U/H$. Moreover $U$  is dense in $Z$ and  contains  the points of $Z$ of codimension $\leq 1$. Furthermore the canonical morphisms $U\to U/H$ is faithfully flat.
\end{theorem}
\proof It has been first stated in \cite{Ray0}, Th\'eor\`eme 1, i) and a proof is contained in \cite{ANA}, Proposition 3.3.1. The last assertion is just \cite{Ray0}, \S 4, Proposition 2. 
\endproof

The conclusion is that, in both cases, $Z \to Z/H$ and $U\to U/H$ are $H$-torsors.

\begin{corollary}\label{corCiVoleva}
Let $S$ be a Dedekind scheme with function field $K$ and $X\to S$ a faithfully flat morphism of finite type. Let moreover $G'$ be an affine and flat $S$-group scheme, $Z$  a faithfully flat $S$-scheme provided with a right $G'$-action $\sigma:Z\times_S G' \to Z$ and $g:Z\to X$ a $G'$-invariant (i.e. $g\circ \sigma=g \circ pr_Z$)  $S$-morphism such that the natural morphism $Z\times_S G' \to Z\times_X Z$ is a closed immersion inducing a $G'_{\eta}$-torsor structure on $Z_{\eta}$ over $X_{\eta}$. Let $U$ be the largest open of $Z$ as in Theorem \ref{teoQuozReno2} such that $U/G'$ is a scheme; then $X':=U/G'$ is faithfully flat and of finite type over $S$ and the natural morphism  $\lambda:X'\to X$ is a model map. In particular $U\to X'$ is a $G'$-torsor  extending the $G'_{\eta}$-torsor $Z_{\eta}\to X_{\eta}$.
\end{corollary}
\proof
By Theorem \ref{teoQuozReno2}, $U$ contains  the points of $Z$ of codimension\footnote{The codimension of a point is defined as the codimension of its closure.} $\leq 1$ so in particular it contains, for all closed points $s\in S$, the generic points of the irreducible components of $Z_s$. As  $U$ is the largest open of $Z$ such that $U/G'$ is a scheme, so in particular it contains $Z_{\eta}$. Thus $U$ is surjective over  $S$. Hence $X'$ is surjective over $S$ too and it is $S$-flat and of finite type because $U$ has the same properties (inherited by $Z$). Thus $X'\to X$ gives rise to the  desired model map. The last assertion is clear.
\endproof

We now recall the definition of N\'eron blow-up.
Unless stated otherwise, from now till the end of section \S \ref{sez:almost} we only consider the following situation:

\begin{notation}\label{notaQF} We denote by $S$ the spectrum of a discrete valuation ring $R$  with uniformising element  $\pi$ and with fraction and residue field respectively denoted by $K$ and $k$. As usual $\eta$ and $s$ will denote  the generic and special point of $S$ respectively. Finally we denote by $X$ a faithfully flat $S$-scheme of finite type.
\end{notation}

According to \cite{BLR}, \S 3.2 Proposition 1 or \cite{ANA}, II, 2.1.2 (A), the following statement holds :

\begin{proposition}\label{propNeronScoppio} Let $S$ be the spectrum of a discrete valuation ring $R$ with uniformising element $\pi$. Let $X$ be a faithfully flat $S$-scheme of finite type and let $C$ be a closed subscheme of the special fiber $X_s$ of $X$ and let $\mathcal{I}$ be the sheaf of ideals of $\mathcal{O}_X$ defining $C$. Let $X'\to X$ be the blow up of $X$ at $C$ and $u:X^C\to X$ denote its restriction to the open subscheme of $X'$ where $ \mathcal{I}\cdot\mathcal{O}_X$ is generated by $\pi$. Then:
\begin{enumerate}
\item $X^C$ is a flat $S$-scheme, $u$ is an affine model map. 
\item For any flat $S$-scheme $Z$ and for any $S$-morphism $v:Z\to X$ such that $v_k$ factors  through $C$, there exists a unique $S$-morphism $v':Z\to X^C$ such that $v=u\circ v'$.
\end{enumerate}
\end{proposition}

\begin{defi}
The morphism $X^C\to X$ (or simply $X^C$) as in Proposition \ref{propNeronScoppio} is called the N\'eron blow up of $X$ at $C$ and property 2 is often referred to as the universal property of the N\'eron blow up. 
\end{defi}

\subsection{Construction of a model}

Given an algebraic $G$-torsor $Y\to X_{\eta}$ we do not know whether or not we can find a torsor over $X$, whose generic fiber is isomorphic to the given one.  In \S \ref{sez:Aim} we have recalled the most important and recent results that partially solve this problem when $G$ is finite; here we suggest a new approach in a much wider context, including the case $G$ of finite type.

\begin{theorem}\label{teoCourbes}
Let $S$ be the spectrum of a discrete valuation ring $R$ with function field $K$. Let $X\to S$ be a faithfully flat morphism of finite type with $X$ a regular and integral scheme of absolute dimension 2 endowed with a section $x\in X(S)$. Let $G$ be a finite $K$-group and $f:Y\to X_{\eta}$ a $G$-torsor pointed in $y\in Y_{x_{\eta}}(K)$, then there exists an integral scheme $X'$ faithfully flat and of finite type over $S$, a model map $\lambda:X'\to X$  and a $G_0$-torsor $f':Y_0\to X'$ extending the given $G$-torsor $Y$ for some quasi-finite and flat $S$-group scheme $G_0$. Moreover $X'$ can be obtained by $X$ after a finite number of N\'eron blow-ups.
\end{theorem}

\begin{theorem}\label{teoDimension}
Let $S$ be the spectrum of a Henselian discrete valuation ring $R$ with function field $K$ and with algebraically closed residue field. Let $X\to S$ be a smooth and surjective morphism with $X$ a connected scheme. Let $G$ be an affine $K$-group scheme of finite type and $f:Y\to X_{\eta}$ a $G$-torsor, then there exists, possibly after extension of scalars, a connected scheme $X'$ smooth and surjective over $S$, a model map $\lambda:X'\to X$  and a $G_0$-torsor $f':Y_0\to X'$ extending the given $G$-torsor $Y$ for some affine finite type and flat $S$-group scheme $G_0$. If moreover $G$ is finite then $G_0$ is quasi-finite and there exists an open sub-scheme $W \subseteq X$ such that $X\backslash W$ has codimension $\geq 3$ in $X$ and such that $X'$ can be obtained from $W$ after a finite number of N\'eron blow-ups.
\end{theorem}

Before proceeding with the proofs of Theorems \ref{teoCourbes} and \ref{teoDimension} we need some preliminary results:

\begin{lemma}\label{lemFinalmente}
Let $S$ be a Dedekind scheme with function field $K$ and $X\to S$ a faithfully flat morphism of finite type with $X$ regular and integral. For any vector bundle $V$  on $X_{\eta}$ there exist a open subscheme $X_1\subseteq X$ containing $X_{\eta}$, where $X\backslash X_1$ has codimension $\geq 3$ in $X$, such that $X_1$ is faithfully flat and of finite type over $S$ and a vector bundle $W$ on $X_1$ such that $W_{|X_{\eta}}\simeq V$. If moreover $dim (X)=2$ then we can choose $X_1=X$
.
\end{lemma}
\proof
Let us denote by $j:X_{\eta}\to X$ the natural open immersion. First of all we observe that there exists a coherent sheaf $\mathcal{F}$ on $X$ such that $j^*(\mathcal{F})\simeq V$ (cf. for instance \cite{H77}, II, ex. 5.15). Then $\mathcal{F}^{\vee\vee}$, i.e. the double dual of  $\mathcal{F}$, is a coherent reflexive sheaf. That $j^*(\mathcal{F}^{\vee\vee})\simeq V$ follows from the well known fact that $j^*(\mathcal{F}^{\vee\vee})\simeq j^*(\mathcal{F})^{\vee\vee}\simeq V$ (see for instance the proof of \cite{Ha}, Proposition 1.8). If $dim(X)=2$ then  by \cite{Ha}, Corollary 1.4 we set $W:=\mathcal{F}^{\vee\vee}$ which is a vector bundle and this is the last assertion. As for the higher dimension case we know, again by  \cite{Ha}, Corollary 1.4, that the subset $C$ of points where $\mathcal{F}^{\vee\vee}$ is not locally free is a closed subset of codimension $\geq 3$. We call $X_1$ the complementary open subset of $C$ in $X$ to which we give the induced scheme structure and we thus know that $X_1$ contains $X_{\eta}$ and has nonempty intersection with $X_s$ for any closed point $s\in S$, whence the first assertion.
\endproof

Let $T$ be any scheme; following \cite{GW} (11.6) we associate to  any locally free sheaf  $V$ of rank $n$ over $T$ the  sheaf $\mathcal{I}som_{\mathcal{O}_T}(\mathcal{O}_T^{\oplus n}, V)$ which is a $GL_{n,T}$-torsor  $\mathcal{I}som_{\mathcal{O}_T}(\mathcal{O}_T^{\oplus n}, V)\to T$, thus obtaining a bijective map between isomorphism classes of locally free sheaves of rank $n$ over $T$ and  isomorphism classes of  $GL_{n,T}$-torsors over $T$.  It is an exercise to prove that this construction base changes \emph{correctly} (i.e. if $i:T'\to T$ is a morphism of schemes then $i^*(\mathcal{I}som_{\mathcal{O}_T}(\mathcal{O}_T^{\oplus n}, V))\simeq \mathcal{I}som_{\mathcal{O}_{T'}}(\mathcal{O}_{T'}^{\oplus n}, i^*(V))$ as $GL_{n,T'}$-torsors).

\proof \emph{of Theorem \ref{teoCourbes}}.  We do the following construction: we take any closed immersion $\rho:G\hookrightarrow GL_{n,K}$ for a suitable $n$ (by  \cite{WW}, \S 3.4). The contracted product $Z:=Y\times^GGL_{n,K}$ via $\rho$
is a $GL_{n,K}$-torsor, then $Z\simeq \mathcal{I}som_{\mathcal{O}_{X_{\eta}}}(\mathcal{O}_{X_{\eta}}^n, V)$ for a suitable vector bundle $V$ on $X_{\eta}$ (for instance one can choose $V:=f_\ast (\mathcal{O} _Y)$) of rank $n$. Let $W$ be a vector bundle on $X$, as in Lemma \ref{lemFinalmente}, whose restriction to $X_{\eta}$ is isomorphic to $V$; let ${Z'}:=\mathcal{I}som_{\mathcal{O}_{X}}(\mathcal{O}_{X}^n, W)$ be the corresponding $GL_{n,S}$-torsor extending $Z$. Let us denote by $\overline{Y}$  and $\overline{G}$ the schematic closures of $Y$ in ${Z'}$ and $G$ in $GL_{n,S}$ respectively. We thus obtain the following diagrams 
\begin{equation}
\label{eqGrande}
\xymatrix{  Y\ar[rr]\ar[dd] \ar@{^{(}->}[rd] & & \overline{Y}\ar[dd]|\hole \ar@{^{(}->}[rd]  &  \\  & Z \ar[rr]\ar[dl] &   & {Z'}\ar[dl] \\  X_{\eta}\ar[rr] \ar[d] & &  X \ar[d] & \\  \eta \ar[rr] & & S &  }\qquad 
\xymatrix{ G\ar[rr]\ar[dd]\ar@{^{(}->}[rd] & & \overline{G}\ar[dd]|\hole \ar@{^{(}->}[rd] & \\  & GL_{n,K}\ar[dl]\ar[rr] &  & GL_{n,S}\ar[dl] \\    \eta\ar[rr] & & S &   }
\end{equation}

In general  neither $\overline{Y}\to X$ nor $\overline{Y}\to S$ will be faithfully flat, however we can modify the embedding $G\hookrightarrow GL_{n,K}$ in order to obtain at least the surjectivity of $\overline{Y}\to S$ (which will be thus faithfully flat). As we will see this will be sufficient to conclude. Let us now pull back the $GL_{n,S}$-torsor $Z'\to X$ over $x:Spec(R)\to X$: we obtain a trivial torsor, whence the existence of a section $z'\in {Z'}_{x}(S)$ whose generic fiber is $z'_{\eta}\in {Z}_{x_{\eta}}(K)$.  Let us now call $z\in {Z}_{x_{\eta}}(K)$ the image of $y$ through the closed embedding $i:Y\hookrightarrow Z$ constructed in diagram (\ref{eqGrande}). In general, unless we are extremely lucky, it will not coincide with $z_{\eta}$, but, as they both live over ${x_{\eta}}$, there exists $g\in GL_{n,K}(K)$ such that $z_{\eta}=z\cdot g$. Let us consider the isomorphism of $X_{\eta}$-schemes $\varphi_g:Z\to Z, z_0\mapsto z_0\cdot g$. The composition $\lambda:=\varphi_g\circ i:Y\hookrightarrow Z$ is a closed immersion sending $y$ to $z_{\eta}$ and it turns out to be a morphism of torsors i.e. commuting with the actions of $G$ and $GL_{n,K}$  if we consider the new embedding $\sigma:G\hookrightarrow GL_{n,K}$ defined as $\sigma:g_0\mapsto g^{-1}\rho(g_0)g$. If we now consider $\overline{Y}'$ and $\overline{G}'$ respectively the  closure of $\lambda: Y\hookrightarrow Z$ in $Z'$ and the closure of $\sigma:G\to  GL_{n,K}$ in $GL_{n,R}$, then we observe that the natural morphism $u:\overline{Y}'\times_S \overline{G}'\to \overline{Y}'\times_{X} \overline{Y}'$ is a closed immersion: it follows from the commutative diagram
$$\xymatrix{{\overline{Y}'}\times_S \overline{G}' \ar[r]^u \ar@{^{(}->}[rd]_i & {\overline{Y}'}\times_{X} {\overline{Y}'}\ar@{^{(}->}[d]^j \\ & {Z'}\times_S GL_{n,S}\simeq {Z'}\times_{X} {Z'}}$$
the fact that both $i:{\overline{Y}'}\times_S \overline{G}'\hookrightarrow {Z'}\times_S GL_{n,S}$ and $j:{\overline{Y}'}\times_{X} {\overline{Y}'}\to {Z'}\times_{X} {Z'}$ are  closed immersions and \cite{ST}, Lemma 28.3.1 (3). Moreover by construction the closure of $y$ in $\overline{Y}'$ has image $z$ in $Z'$ so in particular $\overline{Y}'$ is surjective (and thus faithfully flat) over $S$. According to Theorem \ref{teoQuozReno} we set $X':=\overline{Y}'/\overline{G}'$, $Y_0:=\overline{Y}'$ and $G_0:= \overline{G}'$  in order to conclude. The fact that $X'\to X$ can be obtained as a finite number of N\'eron blow-ups follows from the fact that it is affine and \cite{WWO}, Theorem 1.4. 
\endproof

\proof of the Theorem \ref{teoDimension}
We repeat the first part of the proof of Theorem \ref{teoCourbes} in order to obtain a morphism of torsors $i:Y\hookrightarrow Z$, which is a closed immersion, from the given $G$ torsor to a $GL_{n,K}$-torsor $Z:=\mathcal{I}som_{\mathcal{O}_{X_{\eta}}}(\mathcal{O}_{X_{\eta}}^n, V)$ for some vector bundle $V$ over $X_{\eta}$. Now let $X_1$ be as in Lemma \ref{lemFinalmente} and $W$  a vector bundle on $X_1$ whose restriction to $X_{\eta}$ is isomorphic to $V$; let ${Z_1}:=\mathcal{I}som_{\mathcal{O}_{X_1}}(\mathcal{O}_{X_1}^n, W)$ be the corresponding $GL_{n,S}$-torsor extending $Z$. Let us denote by $\overline{Y}$  and $\overline{G}$ the schematic closures of $Y$ in ${Z_1}$ and $G$ in $GL_{n,S}$ respectively. We thus obtain the following diagrams 
\begin{equation}
\label{eqGrande2}
\xymatrix{  Y\ar[rr]\ar[dd] \ar@{^{(}->}[rd] & & \overline{Y}\ar[dd]|\hole \ar@{^{(}->}[rd]  &  \\  & Z \ar[rr]\ar[dl] &   & {Z_1}\ar[dl] \\  X_{\eta}\ar[rr] \ar[d] & &  X_1 \ar[d] & \\  \eta \ar[rr] & & S &  }\qquad 
\xymatrix{ G\ar[rr]\ar[dd]\ar@{^{(}->}[rd] & & \overline{G}\ar[dd]|\hole \ar@{^{(}->}[rd] & \\  & GL_{n,K}\ar[dl]\ar[rr] &  & GL_{n,S}\ar[dl] \\    \eta\ar[rr] & & S &   }
\end{equation}

As $R$ is Henselian with algebraically closed residue field and as $X_1\to S$ is smooth and surjective there exists a section $x_1\in X_1(S)$. If necessary after a finite extension of scalars $K\hookrightarrow K'$,  $Y_{K'}$ is pointed over $x_{1,\eta}$; we can thus translate the problem to $R'$, the Henselian ring obtained as the integral closure of $R$ in $K'$. In order to ease notations we assume $R=R'$ from now on and we fix a point $y\in Y_{x_{1,\eta}}$. Arguing as in the proof of Theorem \ref{teoCourbes} we can (and we actually do) assume that $i(y)=z_{1,\eta}\in {Z}_{x_{1,\eta}}(K)$ where $z_{1}\in {Z_1}_{x_{1}}(R)$ is a $R$-section of $Z_1\to S$, which always exists as we have seen before, hence $\overline{Y}\to S$ is faithfully flat. As in the proof of Theorem \ref{teoCourbes}, $\bar Y \times _S \bar G \to \bar Y \times _{X_1} \bar Y $ is a closed immersion. Using Corollary \ref{corCiVoleva} there exists a largest open $U$ of $Z$, faithfully flat over $S$, such that $U\to U/\overline{G}$ is a $\overline{G}$-torsor extending the given one. If we set $X':=U/\overline{G}$, $G_0:=\overline{G}$ and $Y_0:=U$ then we obtain the desired result. If $G$ is finite the fact that $U=\overline{Y}$ (we apply to the previous construction Theorem \ref{teoQuozReno} instead of Corollary \ref{corCiVoleva}) implies the last assertion.
\endproof




\begin{thebibliography}{99}

\bibitem{ST} {\sc Stack Project Authors}, \emph{Stacks Project}, version 94a58fd

\bibitem{ANA} {\sc Anantharaman S.},  \emph{Sch\'emas en groupes, espaces homog\`enes et espaces alg\'ebriques sur une base
de dimension 1}. M\'emoires de la S. M. F., tome 33, 5-79 (1973).

\bibitem{Antei4} {\sc  Antei M.}, \emph{Extension of finite solvable torsors over a curve},  Manuscripta Mathematica: Volume 140, Issue 1 (2013), Page 179-194.


\bibitem{Antei2} {\sc  Antei M.}, \emph{On the abelian fundamental group scheme of a family of varieties}, Israel Journal of Mathematics, Volume 186 (2011), 427-446.


\bibitem{AEG} {\sc  Antei M., Emsalem M., Gasbarri C.}, \emph{Sur l'existence du sch\'ema en groupes fondamental}, arXiv:1504.05082 [math.AG].


\bibitem{BLR} {\sc  Bosch S., L\"utkebohmert W., Raynaud M.} \emph{N\'eron models}, Springer Verlag, (1980).



\bibitem{DG} {\sc Demazure M., Gabriel P.}, \emph{Groupes alg\'ebriques}, North-Holland Publ. Co., Amsterdam,
(1970).
%

\bibitem{GAS} {\sc  Gasbarri C.}, \emph{Heights of vector bundles and the fundamental group scheme of a curve}, Duke Math. J. 117, No.2, 287-311 (2003). 
%
\bibitem{GW} {\sc G\"ortz U., Wedhorn T.}, \emph{Algebraic geometry I. Schemes. With examples and exercises.}  Advanced Lectures in Mathematics. Wiesbaden: Vieweg+Teubner

%
%

\bibitem{SGA1} {\sc Grothendieck A.}, \emph{Rev\^etements \'etales et groupe fondamental}, S\'eminaire de g\'eom\'etrie alg\'ebrique du Bois Marie, (1960-61).


\bibitem{DH} {\sc Duong N. G., Hai P. H.}, \emph{Tannakian duality over Dedekind rings and applications}, arXiv:1311.1134v2 [math.AG] 

\bibitem{Ha} {\sc Hartshorne R.}, \emph{Stable Reflexive Sheaves}, Math. Ann. 254, 121-176 (1980)
%
%
\bibitem{H77} {\sc Hartshorne R.}, \emph{Algebraic Geometry}, GTM, Springer Verlag (1977)

\bibitem{MS} {\sc Mehta V. B., Subramanian S.} \emph{The fundamental group scheme of a smooth projective variety over a ring of Witt vectors}, J. Ramanujan Math. Soc. 28A, Spec. Iss., 341-351 (2013).



\bibitem{Mi2} {\sc Milne J. S.}, \emph{Arithmetic duality theorems}, Perspectives in Mathematics, 1. Academic Press, Inc., Boston, MA, (1986).


%
%
\bibitem{Nor2} {\sc Nori M. V.}, \emph{The fundamental group-scheme}, Proc. Indian Acad. Sci. (Math. Sci.), Vol. 91,
Number 2, (1982), p. 73-122.
%
\bibitem{Nor3} {\sc  Nori M. V.}, \emph{The Fundamental Group-Scheme of an Abelian Variety}, Math. Ann. 263, (1983), p. 263-266.
%
%
%

\bibitem{Ray0} {\sc Raynaud M.}, \emph{Passage au quotient par une relation d'\'equivalence plate}, Proceedings of a Conference on Local Fields,
Springer-Verlag (1967), p. 78-85.
%
%


\bibitem{Ray3} {\sc Raynaud M.}, \emph{$p$-groupes et r\'eduction semi-stable des courbes}, The Grothendieck Festschrift, Vol III, Progr. Math., vol. 88, Birkh\"auser, Boston, MA, (1990), p. 179-197.


\bibitem{Ray} {\sc  Raynaud M.}, \emph{Sch\'emas en groupes de type $(p,\ldots,p)$}, Bulletin de la Soci\'et\'e Math\'ematique de France, 102 (1974), p. 241-280. 




\bibitem{Tos} {\sc Tossici D.} \emph{Effective Models and Extension of Torsors over a d.v.r. of Unequal Characteristic}, International Mathematics Research Notices (2008) Vol. 2008 : article ID rnn111, 68 pages (2008).

\bibitem{WW} {\sc  Waterhouse W. C.},
 \emph{Introduction to affine group schemes}, GTM,
Springer-Verlag, (1979).

\bibitem{WWO}  {\sc  Waterhouse W. C., Weisfeiler B.}, \emph{One-dimensional affine group schemes}, Journal of Algebra, 66, 550-568 (1980).

\end{thebibliography}
\end{document}